\documentclass[12pt]{article}
\usepackage{amssymb,latexsym}
\def\ZZ{{\mathbb Z}}

\def\RR{{\mathbb R}}

\newtheorem{formula}{}[section]

\newtheorem{definition}[formula]{\indent Definition}
\newtheorem{corollary}[formula]{\indent Corollary}
\newtheorem{remark}[formula]{\indent Remark}
\newtheorem{lemma}[formula]{\indent Lemma}
\newtheorem{theorem}[formula]{\indent Theorem}

\def\thrm{\begin{theorem}}
\def\thrml#1{\begin{theorem}\label{#1}}
\def\ethrm{\end{theorem}}
\def\rmrk{\begin{remark}}
\def\rmrkl#1{\begin{remark}\label{#1}}
\def\ermrk{\end{remark}}
\def\dfntn{\begin{definition}}
\def\dfntnl#1{\begin{definition}\label{#1}}
\def\edfntn{\end{definition}}
\def\nmrt{\begin{enumerate}}
\def\enmrt{\end{enumerate}}

\def\qtn{\begin{equation}}
\def\qtnl#1{\begin{equation}\label{#1}}
\def\eqtn{\end{equation}}
\def\lmm{\begin{lemma}}
\def\lmml#1{\begin{lemma}\label{#1}}
\def\elmm{\end{lemma}}
\def\crllr{\begin{corollary}}
\def\crllrl#1{\begin{corollary}\label{#1}}
\def\ecrllr{\end{corollary}}

\begin{document}
\title{}
\date{}
\maketitle
\vspace{-0,1cm} \centerline{\bf A CRITERION OF CONTAINMENT FOR TROPICAL HYPERSURFACES}
\vspace{7mm}
\author{
\centerline{Dima Grigoriev}
\vspace{3mm}
\centerline{CNRS, Math\'ematique, Universit\'e de Lille, Villeneuve
d'Ascq, 59655, France} \vspace{1mm} \centerline{e-mail:\
dmitry.grigoryev@univ-lille.fr } \vspace{1mm}
\centerline{URL:\ http://en.wikipedia.org/wiki/Dima\_Grigoriev} }

\begin{abstract}
For tropical $n$-variable polynomials $f, g$ a criterion of containment for tropical hypersurfaces $Trop(f)\subset Trop(g)$ is provided in terms of their Newton polyhedra $N(f), N(g)\subset \RR^{n+1}$. Namely, $Trop(f)\subset Trop(g)$ iff for every vertex $v$ of $N(g)$ there exist a homothety $t\cdot N(f), t>0$ and a parallel shift $s:\RR^{n+1} \to \RR^{n+1}$ such that $v\in s(t\cdot N(f))\subset N(g)$.
\end{abstract}

{\bf keywords}: containment of tropical hypersurfaces, inscribable Newton polyhedra


\section*{Introduction}

Consider a tropical polynomial \cite{MS}
\begin{equation}\label{1}
f=\min_{1\le i\le k} \{M_i\},\ M_i=\sum_{1\le j\le n} a_{i,j}x_j+a_{i,0},\ 0\le a_{i,j}\in \ZZ \cup \{\infty\},\ a_{i,0}\in \RR \cup \{\infty\}.
\end{equation}
The tropical hypersurface $Trop(f)\subset \RR^n$ consists of points $(x_1,\dots,x_n)$ such that the minimum in (\ref{1}) is attained at least at two tropical monomials $M_i, 1\le i\le k$.

For each $1\le i\le k$ consider the ray $\{(a_{i,1},\dots,a_{i,n},a)\ :\ a_{i,0}\le a\in \RR\}\subset \RR^{n+1}$ with the apex at the point $(a_{i,1},\dots,a_{i,n},a_{i,0})$. The convex hull of all these rays for $1\le i\le k$ is Newton polyhedron $N(f)$. Rays of this form we call vertical, and the last coordinate we call vertical. Note that $N(f)$ contains edges (of finite length) and vertical rays. Further, by edges we mean just edges of finite length.

A point  $(x_1,\dots,x_n)\in Trop(f)$ iff a parallel shift $H_x'$ of the hyperplane $H_x=\{(z_1,\dots,z_n,x_1z_1+\cdots+x_nz_n)\ :\ z_1,\dots,z_n\in \RR\}\subset \RR^{n+1}$ has at least two common points (vertices) with $N(f)$, so that $N(f)$ is located in the half-space above $H_x'$ (with respect to the vertical coordinate). In this case $H_x'$ has (at least) a common edge with $N(f)$, and we say that $H_x'$ supports $N(f)$ at $H_x'\cap N(f)$. 

The goal of the paper is to provide for tropical polynomials $f,g$ an explicit criterion of containment $Trop(f)\subset Trop(g)$ in terms of Newton polyhedra $N(f), N(g)$. Note that a criterion of emptiness of a tropical prevariety $Trop(f_1,\dots,f_l)$ is established in \cite{GP} (one can treat this as a tropical weak Nullstellensatz),  further developments one can find in \cite{MR}, \cite{ABG}. The issue of containment of tropical hypersurfaces is a particular case of an open problem of a tropical strong Nullstellensatz, i.e. a criterion of a containment $Trop (f_1,\dots,f_l)\subset Trop(g)$. We mention that in \cite{JM} (which improves \cite{BE}) a strong Nullstellensatz is provided for systems of min-plus equations of the form $f=g$ (in terms of congruences of tropical polynomials). Observe that the family of all tropical prevarieties coincides with the family of all min-plus prevarieties (and both coincide with the family of all finite unions of polyhedra given by linear constraints with rational coefficients \cite{MS}). On the hand, the issue of a strong Nullstellensatz is different for these two types of equations.

\section{Containment of tropical hypersurfaces and inscribable polyhedra}

For a polyhedron $P$ and $0<t\in \RR$ denote by $t\cdot P$ the homothety (with some center) of $P$ with the coefficient $t$.

\begin{definition}\label{inscribe}
For polyhedra $P,Q$ we say that $P$ is inscribed in $Q$ at a point $x$ if $x\in P\subset Q$.

We say that $P\subset \RR^n$ is totally inscribable in $Q$ if for every vertex $v$ of $Q$ an appropriate parallel shift $s:\RR^n \to \RR^n$ of the homothety $s(t\cdot P)$ is inscribed in $Q$ at $v$ for suitable $0<t\in \RR$.
\end{definition}

\begin{theorem}\label{homothety}
For tropical polynomials $f,\ g$ is $n$ variables it holds $Trop(f)\subset Trop(g)$ iff Newton polyhedron $N(f)\subset \RR^{n+1}$ is totally inscribable in $N(g)$.
\end{theorem}

\begin{remark}
Under the conditions of Theorem~\ref{homothety} $s'(t_0\cdot N(f))$ is inscribed in $N(g)$ at an arbitrary chosen point of $N(g)$ (for an appropriate shift $s'$) where $t_0$ is the minimum of $t$ (see Definition~\ref{inscribe}) over all the vertices of $N(g)$ (however, we don't make use of this remark). 
\end{remark}

{\bf Proof of the theorem}. First assume that for every vertex $v$ of $N(g)$ there exists a shift $s$ and $t>0$ such that $s(t\cdot N(f))$ is inscribed in $N(g)$ at $v$. Suppose that $Trop(f)\nsubseteq Trop(g)$, then there exists a hyperplane $\RR^{n+1}\supset H\in  Trop(f) \setminus Trop(g)$. Therefore, a parallel shift of $H$ supports $N(g)$ at some its vertex $v$. By the assumption an appropriate shift $s(t\cdot N(f))$ is inscribed in $N(g)$ at $v$ for suitable $t>0$. This contradicts to that $H\in Trop(f)$ since a parallel shift of $H$ has a single common point $v$ with  $s(t\cdot N(f))$. This proves that $Trop(f)\subset Trop(g)$. \vspace{2mm}

Now conversely, assume that $Trop(f)\subset Trop(g)$. Denote by $p:\RR^{n+1}\twoheadrightarrow \RR^n$ the projection along the last coordinate. Take a vertex $v$ of $N(g)$. Consider a cone $C\subset \RR^{n+1}$ with the apex $v$ being the convex hull of the rays  generated by the edges of $N(g)$ adjacent to $v$ (with the added vertical ray). Then $N(g)\subset C$. Moreover, there exists a ball $B\subset \RR^n$ with the center at $p(v)$ such that $p^{-1}(B)\cap N(g)=p^{-1}(B)\cap C$.

Choose a hyperplane $H\subset \RR^{n+1}$ (not containing a vertical line) such that $H\cap N(g)=\{v\}$, hence $H$ supports $N(g)$ at $v$. Take a vertex $u$ of $N(f)$ for which $H'\cap N(f)=\{u\}$ where $H'$ is a hyperplane parallel to $H$, and $H'$ supports $N(f)$. Observe that $H'\cap N(f)$ is a point since otherwise $H\in Trop(f) \setminus Trop(g)$.

Pick a sufficiently small $t>0$ such that $s(t\cdot N(f))\subset p^{-1}(B)$ where for the shift $s$ holds $s(u_1)=v$, and $u_1$ is the image of $u$ under the homothety (in particular, $v\in s(t\cdot N(f))$). We claim that $s(t\cdot N(f))\subset C$. Indeed, denote by $H_1$ a hyperplane parallel to $H$ and located above $H$. Denote by $L_1,\dots, L_q  \subset \RR^{n+1}$ the rays  with their common apex at $v$ containing edges of $s(t\cdot N(f))$ adjacent with $v$ (with the added vertical ray), and by $C_0 \subset \RR^{n+1}$ the cone generated by $L_1,\dots, L_q$. Then $s(t\cdot N(f))\subset C_0$.

Thus, to justify the claim it suffices to verify that $C_0\subset C$. Suppose the contrary. Denote by $E_1,\dots, E_m$ the rays with their common apex at $v$ containing edges of $N(g)$ adjacent to $v$ (with the added vertical ray), in other words $C$ is the convex hull of $E_1,\dots, E_m$. Denote points $l_i:=L_i\cap H_1, 1\le i\le q,\ e_j:= E_j\cap H_1, 1\le j\le m$. Consider the convex hull $Q\subset H_1$ of the points $l_1,\dots, l_q, e_1,\dots, e_m$. Then a point $l_i$ is one of the vertices of $Q$ for suitable $1\le i\le q$ (according to the supposition). Therefore, there exists a hyperplane $h\subset H_1$ such that $l_i\in h$ and all the points $l_i,\dots, l_{i-1}, l_{i+1},\dots, l_q, e_1,\dots, e_m$ are located in the same of two open half-spaces of $H_1$ separated by $h$. Hence the hyperplane $H_0 \subset \RR^{n+1}$ spanned by $h$ and $v$ belongs to $Trop(g)$, while $H_0 \cap s(t\cdot N(f)) = \{v\}$, i.e. $H_0 \notin Trop(f)$ (observe that $H_0$ does not contain a vertical line since the vertical ray lies in $C\cap C_0$). The obtained contradiction verifies that $C_0\subset C$ and the claim.

Finally, we conclude with
$$s(t\cdot N(f))= s(t\cdot N(f)) \cap p^{-1}(B) \subset C \cap p^{-1}(B) = N(g) \cap p^{-1}(B) \subset N(g).$$ 
\noindent $\Box$ 

\begin{remark}
i) In the proof of Theorem~\ref{homothety} we have chosen a hyperplane $H$ supporting $N(g)$ at a single vertex $v$ in an arbitrary way. On the other hand, a choice of a vertex $u$ of $N(f)$ is subsequently unique (independently of a choice of $H$). Indeed, the space of possible hyperplanes $H$ is connected, and if there were possible to choose another vertex $u_1\neq u$ then for an appropriate choice, $H$ would support $N(f)$ at least at two points, hence $H\in Trop(f) \setminus Trop(g)$.

ii) It would be interesting to provide a criterion of containment for tropical prevarieties $Trop(f_1,\dots,f_k)\subset Trop(g)$. Note that the latter problem is NP-hard \cite{T}, while one can test whether $Trop(f)\subset Trop(g)$ within polynomial complexity (e.g. relying on Theorem~\ref{homothety} and invoking linear programming). 
\end{remark}



\end{document}